\documentclass[english]{amsart}

\usepackage{babel}
\usepackage{amstext}
\usepackage{amsmath}
\usepackage{amsfonts}
\usepackage{latexsym}
\usepackage{ifthen}
\usepackage[all,2cell,dvips]{xy} \UseAllTwocells \SilentMatrices
\xyoption{all}
\newcommand{\codim}{{\rm codim}}

\newcommand{\im}{{\rm im}}

\newtheorem{lemma1}[equation]{}

\newenvironment{lemma}{\begin{lemma1}{\bf Lemma.}}{\end{lemma1}}

\newenvironment{theorem}{\begin{lemma1}{\bf Theorem.}}{\end{lemma1}}
\newenvironment{proposition}{\begin{lemma1}{\bf Proposition.}}{\end{lemma1}}

\newenvironment{definition}{\begin{lemma1}{\bf Definition.}}{\end{lemma1}}

\newenvironment{conjecture}{\begin {lemma1}{\bf Conjecture.}}{\end{lemma1}}

\newcommand{\Z}{\ensuremath{\mathbb{Z}}}
\newcommand{\C}{\ensuremath{\mathbb{C}}}

\newcommand{\PP}{\ensuremath{\mathbb{P}}}

\newcommand{\G}{\ensuremath{\mathbb{G}}}

\newcommand{\Alb}[1]{\ensuremath{\mbox{Alb}(#1)\hspace{0.5ex}}}

\newcommand{\merom}[3]{\ensuremath{#1:#2 \dashrightarrow #3}}

\newcommand{\holom}[3]{\ensuremath{#1:#2  \rightarrow #3}}
\newcommand{\fibre}[2]{\ensuremath{#1^{-1} (#2)}}

\makeatletter
\ifnum\@ptsize=0  \fi \ifnum\@ptsize=2  \fi 

\newcommand\sF{{\mathcal F}}

\newcommand\sO{{\mathcal O}}

\DeclareMathOperator*{\Pic0}{Pic^0}

\DeclareMathOperator*{\nons}{nons}

\newcommand\Ctwo{\ensuremath{C^{(2)}}}

\newcommand{\Chow}[1]{\ensuremath{\mathcal{C}(#1)}}

\newcommand{\Comp}{\ensuremath{\mathcal{K}}}

\newcommand{\barI}{\ensuremath{\bar{I}}}
\newcommand{\Afive}{\ensuremath{\mathcal{A}_5}}
\newcommand{\barAfive}{\ensuremath{\bar{\mathcal{A}_5}}}
\newcommand{\ChowTheta}[1]{\ensuremath{\mathcal{C}_{\frac{\Theta^3}{3!}}(#1)}}
\newcommand{\ChowTwo}{\ensuremath{\mathcal{C}_{\frac{\Theta^2}{2!}}}}
\newcommand{\ChowbarI}{\ensuremath{\mathcal{C}_{\frac{\Theta^3}{3!}}(U/\barI)}}
\newcommand{\ppav}{principally polarised abelian variety}
\newcommand{\ppavs}{principally polarised abelian varieties}

\title{Minimal classes on the intermediate Jacobian of a generic cubic threefold} 
\date{February 7, 2008}

\author{Andreas H\"oring}
\address{Andreas H\"oring, Universit\'e Paris 6, Institut de Math\'ematiques de Jussieu, Equipe de Topologie et G\'eom\'etrie Alg\'ebrique, 175, rue du Chevaleret, 75013 Paris, France}
\email{hoering@math.jussieu.fr}

\begin{document}

\begin{abstract}
Let $X$ be a smooth cubic threefold. We can associate two objects to $X$: 
the intermediate Jacobian $J$ and the Fano surface $F$ parametrising lines on $X$. 
By a theorem of Clemens and Griffiths, the Fano surface can be embedded in the intermediate Jacobian and the cohomology class of its image is minimal. In this paper we show that if $X$ is generic, the Fano surface is the only surface of minimal class in $J$.  
\end{abstract}

\maketitle


\section{Introduction}

Let $X^3 \subset \PP^4$ be a smooth cubic threefold. Its intermediate Jacobian
\[
J(X) := H^{2,1}(X, \C)^* / H_3(X, \Z)
\]
is an abelian variety of dimension five. The intersection product on $H_3(X, \Z)$ induces an irreducible principal polarisation $\Theta$ such that $(J(X), \Theta)$ is not isomorphic to a polarised Jacobian of a curve
\cite[Thm.0.12]{CG72}.
The Fano scheme $F$ parametrising lines contained in $X$ is a smooth surface, and the
Abel-Jacobi map $F \rightarrow J(X)$ is an embedding 
that induces an isomorphism of \ppavs\ $\Alb{F} \simeq J(X)$
\cite[Thm.0.6,0.9]{CG72}. Furthermore the cohomology class of $F \subset J(X)$ is minimal,
that is 
\[
[F] = \frac{\Theta^3}{3!}.
\]
There is only one other known family of 
examples of irreducible principally polarised abelian varieties $(A, \Theta)$ of dimension $n$ such that
for $1 \leq d \leq n-2$, a cohomology class $\frac{\Theta^{n-d}}{(n-d)!}$ can be 
represented by an effective cycle
of dimension $d$: the polarised Jacobians of curves $J(C)$ where the subvarieties 
$W_d(C) \subset J(C)$ have minimal cohomology class. 
In fact on a Jacobian of a curve, these are the only subvarieties having minimal class:
\begin{theorem} \cite[Thm.5.1]{Deb95} \label{theoremdebarre}
Let $C$ be a smooth projective curve of genus $n$, and let $(J(C),\Theta)$ be its Jacobian endowed with the natural principal polarisation. Let $Y \subset J(C)$ be an effective cycle of dimension $1 \leq d \leq n-2$ such that $[Y]=\frac{\Theta^{n-d}}{(n-d)!}$. Then $Y$  is a translate of $W_d(C)$ or $-W_d(C)$.
\end{theorem} 

In this paper we prove a similar statement for the 
intermediate Jacobian $(J(X), \Theta)$ of a {\it generic} smooth cubic threefold. 
In this case the cohomology class 
$\frac{\Theta^{2}}{2!}$ is not represented by an effective cycle \cite[Thm.8.1]{Deb95}.
Since $(J(X), \Theta)$ is not a Jacobian of a curve, the cohomology class 
$\frac{\Theta^{4}}{4!}$ is not represented by an effective cycle by Matsusaka's criterion
\cite{Ma59}. Therefore we are left to identify cycles of class $\frac{\Theta^{3}}{3!}$.

\begin{theorem} \label{maintheorem}
Let $X^3 \subset \PP^4$ be a {\em generic} smooth cubic threefold, and let 
$(J(X), \Theta)$ be its intermediate Jacobian endowed with the natural principal polarisation. 
Let $S \subset J(X)$ be an effective cycle of dimension two that has minimal class, i.e. 
$[S] = \frac{\Theta^{3}}{3!}$.
Then $S$ is a translate of $F$ or $-F$, the Fano surface of $X$.
\end{theorem}

This result provides further evidence for Debarre's conjecture on minimal cohomology classes:

\begin{conjecture} \cite{Deb95} \label{conjecturePP}
Let $(A, \Theta)$ be an irreducible principally polarised abelian variety 
of dimension $n$, and let $Y$ be an 
effective cycle on $A$ of dimension $1 \leq d\leq n-2$.
The following statements are equivalent.
\begin{enumerate}
\item The variety $Y$ has minimal cohomology class, i.e. $[Y] = \frac{\Theta^{n-d}}{(n-d)!}$.
\item $(A, \Theta)$ is the polarised Jacobian of a curve $C$ of genus $n$ and
$Y$ is a translate of $W_d(C)$ or $-W_d(C)$, or $(A, \Theta)$ is the polarised
intermediate Jacobian of a smooth cubic threefold and $Y$ is a translate of $F$ or $-F$. 
\end{enumerate}
\end{conjecture} 

This conjecture holds for \ppavs\ of dimension at most four (the fourfold case is due to Z. Ran \cite[Thm.5]{Ran80}),  
but few things are known about minimal classes in higher dimension.

\medskip

\begin{center}
{\bf
Idea of the proof and structure of the paper.}
\end{center}

Let $S \subset J(X)$ be a subvariety of minimal class, and denote by $S^*$ its smooth locus.
The idea of the proof is to \flqq put a Fano scheme structure on $S^*$\frqq, that is we use the natural inclusion 
\[
\PP(\Omega_{S^*}) \subset \PP(\Omega_{J(X)}|_{S^*}) \simeq S^* \times \PP^4
\]
to interpret $S^*$ as a variety parametrising lines on a threefold $X' \subset \PP^4$.
The hope is of course that we have $X'=X$ up to projective equivalence. 
In fact, the tangent bundle theorem\footnote{We denote by $\PP(E)$ the projectivisation of a vector bundle $E$ in the sense of Grothendieck. This is conform with our main reference \cite{AK77}, but not with Clemens' and Griffiths' seminal paper \cite{CG72}.}\cite[Thm.12.37]{CG72} says that this is true if, up to translation,
$S=\pm F$. The reason why showing $X'=X$ is probably easier than
showing $S=\pm F$ is that we have a characterisation of cubics in $\PP^4$ in terms of the number of lines through a general point.

\begin{theorem} \cite[Thm.1]{La03},\cite{Si30}
Let $X \subset \PP^4$ be an integral hypersurface such that through a general point of $X$ pass only a finite number of lines contained in $X$. If $d$ denotes this number, we have $d \leq 6$. If furthermore $d=6$, then $X$ is a (not necessarily smooth) cubic.
\end{theorem}

Unfortunately it is hard to check the conditions of this theorem in our situation. 
Therefore we will consider a family of
smooth cubic threefolds $X_t$ that degenerate to a cubic $X_0$ that has exactly one ordinary double point: its intermediate Jacobian $J(X_0)$ is an extension of a Jacobian of a curve of genus four by $\C^*$. Using Debarre's Theorem \ref{theoremdebarre} this will allow us to show that a family of cycles of minimal class $S_t$ degenerates to the Fano surface of $X_0$. Once we have this starting point we can deduce the geometrical information we need
to show the following theorem.

\begin{theorem}
\label{theoremdegeneration}
Let $X \rightarrow \Delta$ be a family of three-dimensional cubics over the unit disc
such that for $t \in \Delta^*$ the cubic $X_t$ is smooth and the central fibre
$X_0$ is a cubic with a unique ordinary double point. 
Let $(J,L) \rightarrow \Delta$ be the corresponding family of polarised intermediate Jacobians.
Let $S \rightarrow \Delta$ be a family of irreducible surfaces of minimal class in $J$ such that
$S_0$ is a translate $\pm F_0$, the Abel-Jacobi embedded Fano surface
of $X_0$. 
Then for $t \neq 0$ sufficiently small, $S_t$ is a translate of $\pm F_t$, 
the Abel-Jacobi embedded Fano surface of $X_t$. 
\end{theorem}

This result is our main contribution to the subject. Having defined the objects appearing in the statement,
we will show it in Section \ref{sectionproofdegeneration}.
The moral of the statement is that in order to identify surfaces of minimal class on a generic intermediate Jacobian, 
it is sufficient to do this for a special one.
 
Section \ref{sectionsetup} is devoted to reducing the proof of Theorem \ref{maintheorem}
to Theorem \ref{theoremdegeneration} via fairly standard degeneration arguments.
 
Since the details of our proof are rather technical, 
we should ask if there is a simpler proof using the new techniques introduced by Pareschi and Popa in \cite{PP03,PP06}. In the appendix, 
we show that a surface of minimal class 
in the intermediate Jacobian of a {\em generic} smooth cubic threefold is $M$-regular in their sense
(cf. Proposition \ref{mainproposition}). This is of course a consequence 
of Theorem \ref{maintheorem} and \cite[Thm.]{H07}, but we provide a short and independent proof.

\smallskip

{\bf Acknowledgements.}
I would like to thank Olivier Debarre for many helpful discussions on minimal cohomology classes.

\medskip

\begin{center}
{\bf
Notation}
\end{center}

We work over the complex field $\C$ and denote by $\Delta \subset \C$ the unit disc. 
For standard definitions in algebraic geometry 
we refer to  \cite{Ha77}. 
A variety is an integral complex analytic space of finite dimension, a manifold is a smooth variety.
Points are always closed points.

A fibration is a projective surjective morphism \holom{\phi}{X}{T} with connected fibres
from a complex variety onto a normal complex variety $T$ 
such that $\dim X > \dim T$. 
A fibre is always a fibre in the scheme-theoretic sense, a set-theoretic fibre is the reduction of the fibre.
Fix a complex variety $T$ as a base, 
and let $X \rightarrow T$ be a scheme.  
For $t \in T$ a point, we denote by $X_t$ the fibre over $t$. 
Very often the base $T$ will be the unit disc $\Delta$, in this case $X_0$ will be called the central fibre.

For $(A, \Theta)$ a principally polarised abelian variety, 
we identify $A$ with $\hat{A}=\Pic0(A)$ via the morphism induced by $\Theta$.
We will denote by $P_\xi$ both the points in $\hat{A}=\Pic0(A)$ and the corresponding
numerically trivial line bundle over $A$. We will denote by $\Theta$ the polarisation, i.e. the numerical equivalence
class of an ample divisor, and by $\sO_A(\Theta)$ a divisor representing the polarisation.

For $X$ a projective variety, we denote by $\Chow{X}$ the Chow scheme parametrising effective cycles on $X$ (cf. \cite[Ch.I]{Ko96} for the definition and basic properties).
For simplicity of notation, we will not distinguish between a cycle and its support.

We will say that a property holds for a generic smooth cubic (resp. generic cubic with a unique ordinary double point) if there exists a dense open set in 
the open subset of  $\PP(H^0(\PP^4, \sO_{\PP^4}(3)))$ parametrising smooth cubics (resp. in the subscheme of  $\PP(H^0(\PP^4, \sO_{\PP^4}(3)))$ parametrising 
cubics with a unique ordinary double point) such that the property holds for every variety parametrised by this set.
\newpage

\section{The objects}
\label{sectionobjects}

\subsection{Rank one degenerations of \ppavs}

\begin{definition}\cite{Mu83}
\label{definitionrankone}
A rank one degeneration of a \ppav\ is a pair $(J,L)$, where $J$ is a Gorenstein variety of dimension $n$ and $L$ a
polarisation on $J$ that is constructed as follows:
let $(A, \Theta)$ be a \ppav\ of dimension $n-1$, and let $a \in A$ be a point. 
Recalling that $\Pic0(A)=\mbox{Ext}(A, \C^*)$ as algebraic groups, let 
\[
0 \rightarrow \C^* \rightarrow J^* \rightarrow A \rightarrow 0
\]
be the unique extension associated to the numerically trivial line bundle 
$P_a$. 
Considering $J^*$ as a $\C^*$-bundle over $A$, set 
\[
\tilde{J}:=\PP(\sO_A \oplus P_a)
\] 
for the associated $\PP^1$-bundle, 
and denote by $\holom{\pi}{\tilde{J}}{A}$ the natural map.
We denote by $A_0$ and $A_\infty$ the sections of $\pi$
defined by the projections $\sO_A \oplus P_a \rightarrow \sO_A$
and $\sO_A \oplus P_a \rightarrow P_a$ respectively. 

The variety $J$ is the non-normal variety obtained from $\tilde{J}$ by glueing $A_0$ and $A_\infty$ 
with a translation by $a$. We denote by \holom{\nu}{\tilde{J}}{J} the normalisation morphism.
Set $\tilde{L}:=\pi^* \sO_A(\Theta) \otimes  \sO_{\tilde{J}}(A_\infty)$, 
the polarisation $L$ is represented by the unique line bundle $\sO_J(L)$  such that $\nu^* \sO_J(L) \simeq \tilde{L}$. 
\end{definition}

\begin{definition}
\label{definitionmodulispace}
Let $\Afive$ be the moduli space of \ppavs\ of dimension five. We denote by $\barAfive$ the partial compactification obtained by adding rank one degenerations, 
and by $\partial \Afive$ the locus corresponding to rank one degenerations, i.e. the locus such that
\[
\barAfive = \Afive \sqcup \partial \Afive.
\] 
\end{definition}

{\bf Remark.}
The locus $\partial \Afive$ has a natural fibration $\holom{q}{\partial \Afive}{\mathcal{A}_4}$ 
onto the moduli space of \ppavs\ of dimension four. The fibre over a point $(A, \Theta)$ is the Kummer variety $A/\pm1$ (cf. \cite[Ch.4]{Gru07} for more details on these moduli spaces). 

\subsection{The cubic threefold with a unique ordinary double point}
\label{subsectionsingularcubic1}

Let $X_0$ be a cubic threefold with a unique ordinary double point, and let $F_0:=F(X_0)$ 
be the Fano scheme of lines on $X_0$. 
Then $F_0$ is a non-normal, integral surface that is non-singular except along a  
smooth curve $C_0$ of genus four which is a ordinary double variety of $F_0$ \cite[Thm.7.8]{CG72}.
The curve $C_0$ is a complete intersection of a quadric and a cubic in $\PP^3$, in particular it is not hyperelliptic.
Its normalisation is isomorphic to the symmetric product $\Ctwo_0$ 
and the normalisation morphism $\Ctwo_0 \rightarrow F_0$ consists of identifying two 
disjoint copies $C_0'$ and $C_0''$ of $C_0$ (ibid). 
Note that through a general point of $X_0$ pass only a finite number of lines contained in $X_0$.  

Let $(J(C_0), \Theta_0)$ be the polarised Jacobian of $C_0$, and let
$a$ in $J(C_0) \simeq  \Pic0(\Ctwo_0)$ be 
the point corresponding to  $\sO_{\Ctwo_0}(C_0'-C_0'')$.

\begin{definition} \cite{Co79}
\label{definitionintermediate1}
Let $(J_0, L_0)$ be the rank one degeneration associated to
 $(J(C_0),\Theta_0, a)$.
We call $(J_0, L_0)$ the polarised intermediate Jacobian of $X_0$.\footnote{This notation differs from the one in \cite{CG72}, where the intermediate Jacobian of $X_0$ is only the polarised abelian fourfold $(J(C_0), \Theta_0)$.}
\end{definition}

Let us recall the construction of the embedding $F_0 \rightarrow J_0$.
By \cite[10-12]{Ser58} there exists a morphism 
$\Ctwo_0 \setminus (C_0' \cup C_0'') \rightarrow A$ into an extension of an abelian variety by a torus
that is universal for morphisms from
$\Ctwo_0 \setminus (C_0' \cup C_0'')$ into extensions of abelian varieties by tori.
Using the construction in \cite[11-Thm.1]{Ser58}, one sees 
that this morphism maps to $J_0^*$
where $J_0^*$ is the smooth locus of $J_0$.
Using the notation of Definition \ref{definitionrankone}, we have by construction
\[
J_0^* = \tilde{J}_0 \setminus (A_0 \cup A_\infty),
\]
so we have a commutative diagram
\[
\xymatrix{
& J_0^* \subset \tilde{J}_0 \ar[d]_{\pi}
\\
\Ctwo_0 \setminus (C_0' \cup C_0'') \ar[r] \ar[ru] & J(C_0)
}
\]
Let $\tilde{S}$ be the closure of $\alpha(\Ctwo_0  \setminus (C_0' \cup C_0''))$ in $\tilde{J}_0$.
By \cite[10-01ff]{Ser58} the morphism $\alpha$ is the meromorphic function corresponding to the divisor $\sO_{\Ctwo_0} (C_0'-C_0'')$, so
$A_0 \cap \tilde{S} = C_0''$ and $A_\infty \cap \tilde{S} = C_0'$.
In particular since $C_0$ is not hyperelliptic the birational morphism \holom{\pi|_{\tilde{S}}}{\tilde{S}}{\Ctwo_0} is an isomorphism.
Recall that by definition of $P_a$, we have
\[
P_a \otimes \sO_{\Ctwo_0}  \simeq \sO_{\Ctwo_0} (C_0'-C_0''),
\] 
so $C_0'= a + C_0''$ (seen as subvarieties of the abelian variety $J(C_0)$).
Therefore the restriction of the normalisation morphism \holom{\nu}{\tilde{J}}{J_0} 
to $\tilde{S}$ identifies with the normalisation $\Ctwo_0 \rightarrow F_0$. 

\begin{definition}
\label{definitionfano1}
We will say that  $F_0 \subset J_0$ is the Abel-Jacobi embedded Fano surface of $X_0$.
\end{definition}

\subsection{Fano schemes}
\label{subsectionfanoschemes}

Let $X \rightarrow \Delta$ be a flat family of hypersurfaces of $\PP^4$ over the unit disc, embedded in $\Delta \times \PP^4$. 
Let $\G(1,4)$ be the Grassmannian parametrising projective lines in $\PP^4$.
We denote by $[l]$ 
the point of $\G(1,4)$ parametrising the line $l \subset \PP^4$.
The Fano scheme of $X \rightarrow \Delta$ is defined set-theoretically as
\[
F := F(X/\Delta) := \{ \ (t,[l]) \in \Delta \times \G(1,4) \ | \ l \subset X_t  \} \subset \Delta \times \G(1,4).
\]
By \cite[Thm.3.3.]{AK77} the set $F$ admits a natural projective $\Delta$-scheme structure. Furthermore there exists a universal family, i.e. a scheme 
$$
\Gamma \subset F \times \PP^4 \subset \Delta \times \G(1,4) \times \PP^4
$$ 
endowed with a projection \holom{q}{\Gamma}{F} that makes $\Gamma$ a $\PP^1$-bundle over $F$
such that for every point 
$(t,[l]) \in F$, the fibre \fibre{q}{(t,[l])} is exactly the line 
$l \subset X_t$ parametrised by $(t,[l])$. 
We will denote by 
\[
\holom{\psi}{\Gamma}{\Delta \times \PP^4} 
\]
the restriction of the projection $\Delta \times \G(1,4) \times \PP^4 \rightarrow \Delta \times \PP^4$ to $\Gamma$.
The map $\psi$ is an isomorphism on the fibres of $q$. 
In particular if $x \in X_t \subset \Delta \times \PP^4$ is a point, 
the fibre \fibre{q}{x} 
has finite support if and only if there exist only finitely many lines through contained in $X_t$ and passing through $x$.

The Fano scheme and the family $\Gamma \rightarrow F$ satisfy a universal property:
let $\holom{\phi}{S}{\Delta}$ be a scheme. 
Let $\holom{q_S}{\Gamma_S \subset S \times \PP^4}{\Delta}$ be a $\PP^1$-bundle over $S$ such that for every $s \in S$,
the  fibre \fibre{q_S}{s} is a line $l_s$ that is contained in $X_{\phi(s)} \subset
\PP^4$.
Then there exists a unique $\Delta$-morphism $S \rightarrow F$ such that $\Gamma_S = S \times_F \Gamma$.

\section{Proof of Theorem \ref{theoremdegeneration}.}
\label{sectionproofdegeneration}

{\bf Proof.}
In order to simplify the notation, we will suppose that $S_0$ is a translate of $F_0$, 
the other case is analogous.

Let $J^*$ be the smooth locus of $J \rightarrow \Delta$, i.e. the maximal open subset
such that $\Omega_{J^*/\Delta}$ is locally free. Since the general fibre $J_t$ is an abelian variety
and the central fibre a rank one degeneration, the variety $J^*$ is a group scheme of relative dimension five.
Let $S^*$ be the smooth locus of $S \rightarrow \Delta$.
Note that varieties of minimal cohomology class are reduced, so $S^*$ surjects onto $\Delta$. 
The surjection of relative differential sheaves
\[
\Omega_{J^*/\Delta} \otimes \sO_{S^*} \rightarrow \Omega_{S^*/\Delta} \rightarrow 0
\]
induces an inclusion of projective bundles $\PP(\Omega_{S^*/\Delta}) \hookrightarrow \PP(\Omega_{J^*/\Delta} \otimes \sO_{S^*})$.
For every $t \in \Delta$, the variety $J^*_t$ is an algebraic group, so
$\Omega_{J^*_t}$ is a trivial vector bundle. Up to replacing $\Delta$ 
by a smaller disc, we may therefore suppose that $\PP(\Omega_{J^*/\Delta})$ 
is isomorphic to $\Delta \times J^* \times \PP^4$.
We define the cotangent map by
\begin{equation}
\label{formulacotangentmap}
\phi: \PP(\Omega_{S^*/\Delta}) \hookrightarrow \PP(\Omega_{J^*/\Delta} \otimes \sO_{S^*})
\rightarrow \Delta \times \PP^4,
\end{equation}
where the second arrow is defined by the projection 
$\Delta \times J^* \times \PP^4 \rightarrow \Delta \times \PP^4$.

\smallskip
Since $S^*$ is irreducible, the image $\im(\phi)$ is a constructible, irreducible subset of $\Delta \times \PP^4$ and we denote by $X' \rightarrow \Delta$ 
its closure. 
Since the surfaces $S_t$ are non-degenerate, the image of $\PP(\Omega_{S_t^*})$ has dimension three
\cite[Prop.8.12]{Deb05}.
It follows that $X'$ is an irreducible divisor.
By construction  $X' \rightarrow \Delta$ is surjective, so it is flat
over $\Delta$.
Therefore  $F':=F(X') \rightarrow \Delta$, the Fano scheme of $X'$
exists (cf. Subsection \ref{subsectionfanoschemes}). 
We denote by $\Gamma'$ its universal family together with the natural morphisms $\Gamma' \rightarrow F'$ and
\begin{equation}
\label{equationuniversalmorphism}
\holom{\psi}{\Gamma'}{X'}. 
\end{equation}
By the universal property of the Fano scheme 
there exists a unique rational $\Delta$-morphism
\begin{equation}
\label{formulagaussmap}
\merom{g}{S}{F'} 
\end{equation}
that is defined on $S^*$ such that $\PP(\Omega_{S^*/\Delta}) = S^* \times_{F'} \Gamma'$.
We call $g$ the Gauss map.\footnote{It is not hard to see that these definitions are compatible with the usual definitions on an abelian variety.}

\medskip
{\em We will now study the fibres $X'_t$ of $X' \rightarrow \Delta$ in detail. The goal is to show
that the fibres $X'_t$  identify (up to projective equivalence) to the cubics $X_t$.}
\medskip

Since $\PP(\Omega_{S_t^*})$ is irreducible for every $t \in \Delta$ and
by \cite[Lemma 5.2]{Deb01} the image 
$\phi(\PP(\Omega_{S_t^*}))$ is dense in $X'_t$ for $t \in \Delta$ general, the general fibre
of $X' \rightarrow \Delta$ is an integral hypersurface of $\PP^4$.
Note that a priori it isn't clear if the central fibre $X'_0$ is irreducible or reduced. 
We only know that
it contains the closure of $\phi(\PP(\Omega_{S_0^*}))$ in $\PP^4$. 
By the first statement of Lemma \ref{lemmaevaluationmapcubic} below, this image identifies (up to projective equivalence) to $X_0$.
Therefore we will suppose without loss of generality that 
$X_0$ is an irreducible component of $X'_0$. 
In particular the subvarieties $X'_t \subset \PP^4$ have degree at least three.
{\em We claim that the degree of $X'_t$ is equal to three.} 
\medskip 

{\it Step 1. There are only finitely many lines through a general point of a general fibre $X'_t$.}
We argue by contradiction and suppose that there exist
infinitely many lines through a general point of a general fibre $X'_t$. 
Then the general fibre of $F' \rightarrow \Delta$ has an irreducible component of dimension
at least three, so $F'$ has an irreducible component $F''$ of dimension at least four that
surjects onto $\Delta$. 
Let $\Gamma''$ be the universal family over $F''$, then $\Gamma''$ 
is irreducible of dimension at least five.

Consider now the map \holom{\psi|_{\Gamma''}}{\Gamma''}{X'}, obtained as the restriction of the morphism \ref{equationuniversalmorphism} to $\Gamma''$; it is a $\Delta$-morphism that is dominant by construction, 
thus surjective by the properness of $\Gamma''$ and the irreducibility of $X'$. 
Furthermore since $X'$ has dimension four, the fibres of $\psi|_{\Gamma''}$ have dimension at least one.
In particular if we choose a general point $x \in X_0$, the fibre is not finite.
As noted in Subsection \ref{subsectionfanoschemes}, this implies that 
through a general point of $X_0 \subset X_0'$ there exists an infinity of lines contained in $X_0$, 
a contradiction.

{\it Step 2. Proof of the claim.}
Note that since the cubic $X_0$ is an irreducible component of $X'_0$, one of the irreducible components
of $F'_0$ is the Fano surface of $X_0$.

We consider the Gauss map \merom{g}{S}{F'}, and denote by $F''$ the reduction of the irreducible
component of $F'$ that contains $g(S^*)$. 
Let $\Gamma''$ be the universal family over $F''$, and denote by \holom{\psi|_{\Gamma''}}{\Gamma''}{X'} the restriction of the morphism \ref{equationuniversalmorphism} to $\Gamma''$.
Since the lines parametrised by $S^*$ dominate the general fibres of $X'$, the map $\psi|_{\Gamma''}$ is surjective.
Since there are only finitely many lines through a general point of a general fibre $X'_t$,
the morphism $p$ is generically finite. 

Step 1 allows us to apply Landsberg's theorem \cite[Thm.1]{La03}: there are at most six lines through a general point of a general fibre. Thus the morphism $\psi|_{\Gamma''}$ has degree at most six.
On the other hand the image of $g(S^*_0)$ is dense in the Fano surface of $X_0$, so  
we have a set-theoretical inclusion $F_0 \subset F''$.
Denote by $\Gamma_0 \subset \Gamma''$ the universal family over $F_0$.
Since $X_0$ has only an isolated singularity, we know by \cite[Prop.1.7]{AK77} that there exist six
lines through a general point of $X_0$. Hence for $x \in X_0$ general, the fibre of the restricted morphism
\[
\holom{\psi|_{\Gamma_0}}{\Gamma_0 \hookrightarrow \Gamma''}{X'}
\]
consists of six points, so $\psi|_{\Gamma''}$ has degree at least six.

Hence there are exactly six lines through a general point of a general fibre $X'_t$.  By \cite{Si30}  this implies that the threefold $X'_t$ has degree three. This proves the claim.

{\it Step 3. Structure of the fibres $X'_t$.} 
The general fibre has degree three, so by flatness all the fibres have degree three. Since the cubic $X_0$
is an irreducible component of the central
fibre $X_0'$, the fibre $X_0'$ is integral and identifies to $X_0$.
Up to replacing $\Delta$ by a smaller disc, we
can suppose without loss of generality that all the fibres $X_t'$ are integral cubic hypersurfaces. 
Since the central fibre has exactly one ordinary double point, its Milnor number equals one \cite[I, Thm.2.46]{GLS06}.
Thus the upper semicontinuity of the Milnor number \cite[I, Thm.2.6.]{GLS06}
shows that the general fibre is nonsingular or has exactly one ordinary double point.  

Let us now compare the families $S \rightarrow \Delta$ and $F' \rightarrow \Delta$.
The central fibres $S_0$ and $F_0'$ 
are both isomorphic to the Fano surface $F_0$ of $X_0$. 
By the second statement of Lemma \ref{lemmaevaluationmapcubic}, the rational map $g|_{S_0}$ extends to an isomorphism $S_0 \simeq F_0'$.
Since the indeterminacy locus of $g$ is closed over $\Delta$, this shows that, up to replacing $\Delta$ by a smaller disc, 
$g$ extends to a finite morphism defined on $S$. 
In order to see that the degree of $g$ is one, we have to show that $g$ does not ramify in a general point of $S_0$.

The central fibre of $S \rightarrow \Delta$ is reduced, so if $s \in S_0$ is a smooth point,
the kernel of the tangent map
\[
\holom{Tg_s}{T_{S,s}}{T_{F',g(s)}} 
\]
is contained in $T_{S_0,s}$. Yet $Tg_s|_{T_{S_0,s}}$  is injective, since the restriction $g|_{S_0}$ is a local isomorphism. Therefore the tangent map is injective, so $g$ has degree one. 
In particular it is an isomorphism onto the smooth locus of $F'$.

We will now argue by contradiction and suppose that the general fibre of $X' \rightarrow \Delta$ has an ordinary double point. Then for $t \neq 0$, the Fano schemes $F'_t$ are non-normal, integral surfaces and their  
normalisation is isomorphic to $C_t^{(2)}$, with 
$C_t$ a smooth curve of genus four (cf. Subsection \ref{subsectionsingularcubic1}). 
The map $g|_{S_t}$ is an isomorphism onto the nonsingular locus $(F'_t)_{\nons} \subset C_t^2$, so we can define
the inverse rational map $C_t^{(2)} \dashrightarrow S_t$. 
Since $C_t^{(2)}$ is smooth and $S_t$ contained in the abelian variety $J(X_t)$, the rational map extends to a surjective morphism 
$C_t^{(2)} \rightarrow S_t$ \cite[Thm.4.9.4]{BL04}. 
By the universal property of the Albanese map this morphism factors through the Albanese map $C_t^{(2)} \rightarrow \Alb{C_t^{(2)}}$. It follows that $S_t$ is contained in the image of
the morphism  $\Alb{C_t^{(2)}} \rightarrow J(X_t)$ which has dimension at most four. This contradicts the non-degeneracy of $S_t$. Hence the general fibre of $X' \rightarrow \Delta$ is smooth.

{\it Step 4. Structure of the surfaces $S_t$.}
Since for $t \neq 0$, the cubic $X'_t$ is smooth, the Fano surface $F'_t$ is smooth. 
Since the Gauss map \holom{g}{S}{F'} is an isomorphism onto the smooth locus of $F'$, 
the restriction \holom{g|_{S_t}}{S_t}{F_t'} is an isomorphism.
Hence $S_t \simeq F_t'$, and the morphism $F'_t \rightarrow S_t \subset J(X_t)$ factors through the Albanese map of $F'_t$.
The induced map $\Alb{F'_t} \rightarrow J(X_t)$ is an isomorphism of \ppavs, 
so by the Torelli theorem for smooth cubic threefolds \cite[Thm.0.11]{CG72} we have $X_t=X'_t$ (up to projective equivalence). Since the Albanese map $F'_t \hookrightarrow \Alb{F'_t}$
identifies to an Abel-Jacobi map $F'_t \hookrightarrow J(X'_t)$, the surface $S_t$ is Abel-Jacobi embedded.  
$\square$

\begin{lemma} \cite{CG72, AK77}. 
\label{lemmaevaluationmapcubic}
Let $J_0$ be the intermediate Jacobian of a cubic threefold $X_0$
with exactly one ordinary double point, and let $F_0 \subset J_0$ 
be the Abel-Jacobi embedded Fano surface of $X_0$.  
\begin{enumerate}
\item Let \holom{\phi_0}{\PP(\Omega_{(F_0)_{\nons}})}{\PP^4} be the cotangent map (\ref{formulacotangentmap}), and
let $X_0'$ be the closure of its image. Then $X_0'=X_0$ (up to projective equivalence).
\item Let \merom{g_0}{F_0}{F_0':=F(X_0')} be the Gauss map (\ref{formulagaussmap}).
Then $g_0$ extends to an isomorphism $F_0 \simeq F_0'$. 
\end{enumerate}
\end{lemma}

{\bf Proof.}
Statement 1) is a mere reformulation of \cite[Thm.1.10,ii)]{AK77}. 

To show statement 2), note that by the first statement 
the varieties $F_0$ and $F_0'$ are isomorphic as abstract varieties, so the question is only if
the birational map $g_0$ extends to an isomorphism. Since $g_0$ is defined by sending a point in $F_0$ to the point
parametrising the corresponding line in $X_0'$, this is somewhat obvious. More formally, we can argue
that if \merom{\tilde{g}_0}{\tilde{F}_0}{\tilde{F}'_0} is the rational map induced by $g_0$ on the normalisations,
then $\tilde{F_0} \simeq \Ctwo_0$ is smooth and $\tilde{F}'_0 \simeq \Ctwo_0$ contains no rational curves
(recall that $C_0$ is not hyperelliptic). Therefore $\tilde{g}_0$ extends to a morphism on $\tilde{F}_0$  by 
\cite[Cor. 1.44]{Deb01}, hence $g_0$ extends to a morphism on $F_0$. 
$\square$ 

\newpage
\section{The degeneration argument}
\label{sectionsetup}

\subsection{Technical setup}
\label{subsectionchow}

This section is mainly concerned with the setup of our degeneration argument. 
The material we use is well-known to experts and has been used in the same way by various authors (e.g. \cite{BC77,Deb95}). 
For the convenience of the reader we nevertheless write down the details.  
\smallskip

Let $I \subset \Afive$ be the locus of polarised intermediate Jacobians of smooth cubic threefolds in $\PP^4$, the scheme $I$ is irreducible of dimension ten.
Denote by $\barI$ the closure of $I$ in $\barAfive$. Then
$\barI$ contains the locus of Jacobians of hyperelliptic curves of genus five \cite[Thm.1]{Cas06}.
Let $q$ be the natural fibration on $\partial \Afive$ (cf. Definition \ref{definitionmodulispace}). 
The image $q(\partial \Afive \cap \barI)$ is the closure of the locus of polarised Jacobians of curves of genus four \cite[p.333]{Deb95}.

Up to replacing $\barI$ by a finite covering, we can suppose that there exists a versal family
$U \rightarrow \barI$ endowed with a polarisation $L \rightarrow U$ such that for $t \in \barI$, the couple 
$(U_t, L|_{U_t})$ is isomorphic to the polarised variety parametrised by the point $t$.
The relative Chow functor parametrising cycles of dimension two and fixed degree $\frac{5!}{3!}$ with respect to this polarisation is representable
by a scheme that is projective over $\barI$ \cite[I, Thm3.21]{Ko96}. 
Let $\ChowbarI$ be the subscheme that parametrises the two-dimensional cycles of minimal class. 
By \cite[Ch.6]{Deb95},
the natural morphism  $\holom{p}{\ChowbarI}{\barI}$ is projective and surjective.

Let $\Comp \subset \ChowbarI$ be the reduction of an irreducible component such that the natural morphism $\holom{p:=p|_{\Comp}}{\Comp}{\barI}$ is surjective. If $z \in \fibre{p}{t}$ is a general point in a general fibre $\fibre{p}{t}$ and we cut $\Comp$ with sufficiently general very ample divisors, we get a subvariety
$Z$ passing through $z$ such that $\holom{p|_Z}{Z}{\barI}$ is generically finite.
Up to making a base change $Z \rightarrow \barI$, 
we may suppose without loss of generality that $p|_Z$ is birational. 
Let now $\holom{\theta}{\Delta}{\barI}$ be a finite morphism such that $0$ maps to a certain
point $t_0 \in \barI$.  For a generic choice of $\theta$, the image $\theta(\Delta)$ meets the locus
on which  $\holom{p|_Z}{Z}{\barI}$ is an isomorphism.
So, up to replacing $\Delta$ by a smaller disc, we have a section 
$s: \Delta^* \rightarrow \Comp \times_{\barI} \Delta$.  We extend the corresponding family
$S^* \rightarrow \Delta^*$ of surfaces of minimal class to a family of surfaces $S$ parametrised by $\Delta$. 

{\it Throughout the rest of the paper, we will suppose (without loss of generality) that 
given $\Comp$ and a morphism $\holom{\theta}{\Delta}{\barI}$ as above, we can consider
a corresponding family of surfaces. \label{convention}}

\smallskip
We start our study the Chow scheme with a technical lemma.

\begin{lemma}
\label{lemmatechnical}
Let $X^3 \subset \PP^4$ be a {\em generic} smooth cubic threefold, and 
let $(J(X), \Theta)$ be its polarised intermediate Jacobian.
Let $\mathcal{C}_{\frac{\Theta^3}{3!}}$ be the subscheme of $\Chow{J(X)}$ parametrising cycles of cohomology class $\frac{\Theta^{3}}{3!}$.
Then the connected components of  $\mathcal{C}_{\frac{\Theta^3}{3!}}$ are irreducible of dimension five. If $S$ and $S'$ are two cycles parametrised by the same component, they are translates of each other. 
\end{lemma}

{\bf Proof.}
Since $\ChowbarI$ has {\em finitely many} irreducible components, understanding
cycles of minimal class on the intermediate Jacobian of a {\em generic} smooth cubic threefold is equivalent to understanding the cycles parametrised by the irreducible components of $\ChowbarI$ that surject onto $\barI$.

We will need the following construction: 
let $t \in \barI$ be a point parametrising an irreducible \ppav\ $(A, \Theta)$. 
Fix a point $z_0$ in the fibre $\ChowbarI_t$, and let $S_{z_0}$ be the cycle of minimal
class parametrised by the point. The translates of $S_{z_0}$
form a family of cycles of minimal class parametrised by $A$, and we denote by $Z(z_0) \subset \ChowbarI_t$
the corresponding irreducible subscheme of the Chow scheme. 
The scheme $Z(z_0)$ has obviously dimension at most five, and in fact equal to five:
otherwise there would exist a subvariety $D \subset A$ of positive dimension such that for all $a \in D$, we have
\[ 
S_{z_0} + a = S_{z_0}.
\]
In particular $S_{z_0}$ would be geometrically degenerate \cite[Ch.VIII, Cor.11]{Deb05} which is impossible for effective cycles of minimal class \cite[Ch.1]{Deb95}. 

{\it Claim. The irreducible components of a general $\ChowbarI_t$ are isomorphic to some $Z(z_0)$.}
Let $t \in \barI$ be a point corresponding to a Jacobian $(J(C), \Theta)$ 
of a hyperelliptic curve $C$.
The Chow scheme $\ChowTheta{J(C)}$ parametrising cycles of cohomology class $\frac{\Theta^{3}}{3!}$ has exactly one irreducible component: since $C$ is hyperelliptic, Debarre's Theorem \ref{theoremdebarre} implies that the cycles of this class are translates of $W_2(C)=-W_2(C)$.
In particular $\ChowTheta{J(C)}$ has dimension five.
The upper semicontinuity of the fibre dimension applied to all the irreducible components of $\ChowbarI$ that surject onto $\barI$ shows that a general $\ChowbarI_t$ has dimension at most five.
By the construction above, we see that for $t \in \barI$ general, 
every point $z_0$ in the fibre $\ChowbarI_t$ is contained in $Z(z_0) \subset \ChowbarI_t$,
which has dimension five. A straightforward argument implies that the irreducible components of $\ChowbarI_t$ are of the form $Z(z_0)$.  This shows the claim.

The translation induces an equivalence relation on the cycles, so two 
subschemes $Z(z_0)$ and $Z(z_1)$ are either disjoint or identical. Therefore the claim implies that two cycles parametrised by the same components are translates of each other. $\square$

\begin{lemma} \label{lemmaintegral}
In the situation of the preceding lemma,
an effective cycle of class $\frac{\Theta^{3}}{3!}$ is integral.
\end{lemma}

{\bf Proof.}
We argue by contradiction and suppose that there exists an irreducible component $\Comp \subset \ChowbarI$ that surjects onto $\barI$, and such that for $t \in \barI$ generic, the fibre $\Comp_t$ has an irreducible component that parametrises cycles that are not integral.
Let now $\holom{\theta}{\Delta}{\barI}$ be a finite morphism such that $0$ maps to a point corresponding to a polarised Jacobian of a hyperelliptic curve $C$. 
For a generic choice of $\theta$, we consider
a corresponding family $S \rightarrow \Delta$ of surfaces of minimal class (cf. page \pageref{convention}). 
The central fibre $S_0$ is of minimal class in a hyperelliptic Jacobian, 
so it is a translate of the variety $W_2(C)$ by Debarre's Theorem \ref{theoremdebarre}.
The central fibre is integral and the base is irreducible, so $S$ is irreducible.
Since integrality is an open property in irreducible families, 
we see that a general cycle parametrised by $\Delta$ is integral, a contradiction. 

This shows the claim for a cycle
parametrised by a generic point of a generic $\Comp_t$. By Lemma \ref{lemmatechnical},
any point in a generic $\Comp_t$ corresponds to the translate of a cycle parametrised by a generic point, so the full claim follows.
$\square$

The results obtained so far 
were based on degenerating a family of polarised intermediate Jacobians $(J(X_t), \Theta_t)$ into a polarised Jacobian of a hyperelliptic curve
of genus five $(J(C_0), \Theta_0)$. For the geometric arguments that we want to make in the following, it is not possible to consider this type of degeneration: 
as shown by Collino \cite{Co82}, the degeneration into the hyperelliptic Jacobian corresponds
to degenerating the smooth cubics $X_t$ into a cubic $X_0$ that is singular along a rational normal curve of degree four. The Fano scheme of such a cubic is a union of two
copies of $\PP^2$, in particular we can't \lq see\rq\ the Fano scheme of $X_0$ in the hyperelliptic Jacobian.\footnote{We have seen above that the Abel-Jacobi embedded Fano surfaces $F_t$ degenerate into $W_2(C_0)$. For a geometric explanation of this degeneration, cf. \cite[Prop.2.1]{Co82}.} For this reason, we will now degenerate the smooth cubics  
$X_t$ into a cubic $X_0$ that has a unique ordinary double point. The intermediate Jacobian
$J(X_0)$ is only semi-abelian, but on the level of the Fano schemes, the situation will be much better.

\begin{lemma}
\label{lemmaidentifyF0}
Let $X \rightarrow \Delta$ be a family of three-dimensional cubics over the unit disc
such that for $t \in \Delta^*$ the fibre $X_t$ is a general smooth cubic and the central fibre
$X_0$ is a general cubic with a unique ordinary double point. 
Let $(J,L) \rightarrow \Delta$ be the corresponding family of  polarised intermediate Jacobians.  
Denote by $\holom{p}{\mathcal{C}_{\frac{\Theta^3}{3!}}(J/\Delta)}{\Delta}$ the relative Chow scheme parametrising surfaces of minimal class, and fix an irreducible component  $\Comp \subset \mathcal{C}_{\frac{\Theta^3}{3!}}(J/\Delta)$ that surjects onto $\Delta$. 

Let $S \rightarrow \Delta$ be a family of surfaces of minimal class in $J$ such that $S_0$ corresponds
to a generic point of $\Comp_0:=\fibre{p}{0} \cap \Comp$. 
Then the surface $S_0$ is a translate of $\pm F_0$, where $F_0 \subset J_0$ is the Abel-Jacobi embedded Fano surface of $X_0$. 
\end{lemma}

{\bf Proof.} 
Let \holom{\nu}{\tilde{J}_0}{J_0} be the normalisation map.
By the Definitions \ref{definitionrankone} and \ref{definitionintermediate1}, 
we have
\[
\tilde{J}_0 \simeq \PP(\sO_{J(C_0)} \oplus P_a)
\] 
where $(J(C_0),\Theta_0)$ is a Jacobian of a non-hyperelliptic curve of genus four and $P_a$ is a numerically trivial line bundle over $J(C_0)$. We denote by \holom{\pi}{\tilde{J}_0}{J(C_0)} the projection map, and by 
$A_0$ and $A_\infty$ the $\pi$-sections defined in \ref{definitionrankone}.
Let $\tilde{S}_0 \subset \tilde{J}_0$ be the closure of the inverse image of $S_0$ under the bijection
\[
\holom{\nu|_{(\tilde{J}_0 \setminus A_\infty)}}{(\tilde{J}_0 \setminus A_\infty)}{J_0}.
\]
The proof of \cite[Thm.8.1, Lemma 8.2]{Deb95} shows that there are two cases:
\begin{enumerate}
\item either $\tilde{S}_0$ is irreducible. 
Then the morphism \holom{\pi|_{\tilde{S}_0}}{\tilde{S}_0}{\pi(\tilde{S}_0)} is an isomorphism  onto a translate of $\pm \Ctwo_0$.
Furthermore $\tilde{S}_0$ meets $A_0$ (resp. $A_\infty$) along a translate of $\pm C_0$ (resp. $\pm C_0-a$) and $S_0$ is obtained from $\tilde{S}_0$ by glueing these intersections.
\item or $\tilde{S}_0$ has two irreducible components $S_1$ and $S_2$ where 
\begin{itemize}
\item $S_1 \subset A_0$ and  
\holom{\pi|_{S_1}}{S_1}{\pi(S_1)} is an isomorphism onto a translate of $\pm \Ctwo_0$.
\item $S_2 = \fibre{\pi}{V}$ where $V$ is a translate of $\pm C_0$ such that $V \subset \pi(S_1)$ 
and $V-a \subset \pi(S_1)$.
\end{itemize}
\end{enumerate}

\smallskip
{\em We claim that if $S_0$ corresponds to a generic point of $\Comp_0$, we are in the first case.}

\smallskip
Assume this for the time being. 
Then the surface $\tilde{S}_0$ is isomorphic to $\Ctwo_0$. 
Furthermore $S_0$ is obtained from $\tilde{S}_0$ by 
glueing $\tilde{S}_0 \cap A_0$  and $\tilde{S}_0 \cap A_\infty$ 
via the translation $x \rightarrow x-a$.
This shows that $S_0$ and $F_0$ are isomorphic and
the inclusion $S_0 \rightarrow J_0$ is the Abel-Jacobi embedding (cf. Definition \ref{definitionfano1}).
By the universal property $S_0$ is a translate of $\pm F_0$.

\smallskip
{\em Proof of the claim.}
Let $\tilde{Z}$ be the subset the Chow scheme of $\tilde{J}_0$
corresponding to the cycles $\tilde{S}_0$ described above. We have 
\[
\Comp_0 \subset \nu_* \tilde{Z},
\]
where $\nu_*$ is the morphism on the Chow schemes induced by $\holom{\nu}{\tilde{J}_0}{J_0}$.
By Lemma \ref{lemmatechnical} and upper semicontinuity of the fibre dimension, the irreducible components
of $\Comp_0$ have dimension at least five. We are done if we show that the irreducible components
of $\tilde{Z}$ have dimension five and a general point of $\tilde{Z}$ corresponds to an irreducible cycle.

The morphism $\holom{\pi}{\tilde{J}_0}{J(C_0)}$
induces a morphism \holom{\pi_*}{\tilde{Z}}{\ChowTwo} to the scheme $\ChowTwo$ parametrising the translates of $\pm \Ctwo_0$. 
This scheme has two irreducible components of dimension four. Let $c \in \ChowTwo$ 
be a generic point. The fibre $\fibre{\pi_*}{c}$
parametrises the cycles $\tilde{S}_0$ that map onto a fixed translate $S'$ of $\pm \Ctwo_0$. 
The surface $\tilde{S}_0$ is a divisor in the smooth threefold 
\[
\fibre{\pi}{S'} \simeq \PP(\sO_{S'} \oplus P_a|_{S'})
\]
and the dimension of the fibre $\fibre{\pi_*}{c}$ at the point corresponding to $\tilde{S}_0$ 
equals the dimension of the linear system $|\sO_{\fibre{\pi}{S'}}(\tilde{S}_0)|$.
Consider the exact sequence
\[
0 \rightarrow \sO_{\fibre{\pi}{S'}}(\tilde{S}_0-A_\infty) \rightarrow \sO_{\fibre{\pi}{S'}}(\tilde{S}_0)
\rightarrow \sO_{\fibre{\pi}{S'} \cap A_\infty}(\tilde{S}_0) \rightarrow 0.
\]
The intersection $A_\infty \cap \tilde{S}_0$ is a translate of $\pm C_0$. 
In order to simplify the notation we will
assume that the intersection is a translate of $C_0$. 
Since $\fibre{\pi}{S'} \cap A_\infty  \simeq \Ctwo_0$, 
we see that $h^0(\fibre{\pi}{S'} \cap A_\infty, \sO_{\fibre{\pi}{S'} \cap A_\infty}(\tilde{S}_0))= h^0(\Ctwo_0, \sO_{\Ctwo_0}(C_0))=1$.
Since $\sO_{\tilde{J}_0}(A_0-A_\infty) \simeq \pi^* P_a$, the description of $\tilde{S}_0$ yields
\[
\sO_{\fibre{\pi}{S'}}(\tilde{S}_0-A_\infty) \simeq \pi^* (\sO_{S'}(C_0) \otimes P_a|_{S'}),
\]
so again
\[
h^0 (\fibre{\pi}{S'},\sO_{\fibre{\pi}{S'}}(\tilde{S}_0-A_\infty))
= h^0 (S', \sO_{S'}(C_0) \otimes P_a|_{S'})=1.
\]
The long exact cohomology sequence shows that
the linear system $|\sO_{\fibre{\pi}{S'}}(\tilde{S}_0)|$ has dimension equal to one, so
the components of $\tilde{Z}$ have dimension five. Furthermore
the fixed part of the linear system  $|\sO_{\fibre{\pi}{S'}}(\tilde{S}_0)|$ is empty, so a general member is irreducible. $\square$

\subsection{Proof of Theorem \ref{maintheorem}}
\label{subsectionproof}

Consider the map $\holom{p}{\ChowbarI}{\barI}$ introduced in Subsection \ref{subsectionchow}. 
We argue by contradiction and suppose that there exists an irreducible component 
$\Comp \subset \ChowbarI$ that surjects onto $\barI$, 
and such that for $t \in \barI$ generic, the fibre $\Comp_t$ has an irreducible component 
such that the cycles parametrised by this component are not translates of the Fano surface.
Choose a generic finite morphism
$\holom{\theta}{\Delta}{\barI}$ such that $0$ 
maps to a point corresponding to a polarised intermediate Jacobian of a cubic with a unique ordinary double point, 
and let $(J,L) \rightarrow \Delta$ be the corresponding family of polarised varieties. 
Take a family $S \rightarrow \Delta$ of surfaces of minimal class in $J$ such that
$S_0$ corresponds to a generic point of $\Comp \times_{\barI} \Delta$ (cf. page \pageref{convention}).
By Lemma \ref{lemmaintegral}, the general surface $S_t$ is integral.
Moreover the surface $S_0$ is a 
translate of the Fano surface $\pm F_0$ by Lemma \ref{lemmaidentifyF0}.
By Theorem \ref{theoremdegeneration} the general surface $S_t$ is a translate of the Fano surface $\pm F_t$, a contradiction. $\square$

\section{Appendix: M-regularity}
\label{sectionmregularity}

\begin{definition} \cite{PP06},\cite[Lemma 3.8]{PP06b}
\label{definitiongvsheaf}
Let $(A, \Theta)$ be an irreducible \ppav\ 
of dimension $n$, and let $\sF$ be a coherent sheaf
on $A$. For all $1 \leq i \leq n$, we denote by
\[
V^i := \{ \xi \in \Pic0(A) \ | \ h^i(A, \sF \otimes P_\xi)>0 \}
\]
the $i$-th cohomological support locus of $\sF$. 
We say that $\sF$ is $M$-regular if
\[
\codim V^i > i 
\]
for all $1 \leq i \leq n$.
\end{definition}

The following proposition could be a starting point for an alternative proof of Theorem \ref{maintheorem}.

\begin{proposition} \label{mainproposition}
Let $X^3 \subset \PP^4$ be a {\em generic} smooth cubic threefold, and 
let $(J(X), \Theta)$ be its polarised intermediate Jacobian.
Let $S \subset J(X)$ be an effective cycle of dimension two that has minimal cohomology class, i.e. 
$[S] = \frac{\Theta^{3}}{3!}$.
The surface $S$ is integral, the twisted structure sheaf $\sO_S(\Theta)$ is $M$-regular, and 
$h^0 (S, \sO_S(\Theta) \otimes P_\xi) = 1$ for $P_\xi \in \Pic0(A)$ general. 
\end{proposition}

{\bf Proof of Proposition \ref{mainproposition}.}
The integrality of the cycles has already been proven in Lemma \ref{lemmaintegral}.

Consider the map $\holom{p}{\ChowbarI}{\barI}$ introduced in Subsection \ref{subsectionchow}. 
We argue by contradiction and suppose that 
there exists an irreducible component $\Comp \subset \ChowbarI$ that surjects onto $\barI$, and such that for $t \in \barI$ generic, the fibre $\Comp_t$ has an irreducible component 
such that the claim is false for cycles parametrised by this component.
We proceed as in the proof of Lemma \ref{lemmaintegral}: choose a generic finite morphism
$\holom{\theta}{\Delta}{\barI}$ such that $0$ 
maps to a point corresponding to a polarised Jacobian $(J(C_0), \Theta_0)$ 
of a hyperelliptic curve, and let $(J,L) \rightarrow \Delta$ be the corresponding family
of \ppavs. 
By our assumption we can take a family $S \rightarrow \Delta$ of surfaces of minimal class in $J$ such that
for $t \neq 0$ the claim does not hold for $S_t$ (cf. page \pageref{convention}). 

The central fibre $S_0$ is of minimal class in a hyperelliptic Jacobian, 
so it is a translate of the variety $W_2(C_0)$ by Debarre's Theorem \ref{theoremdebarre}.
Therefore by \cite[Prop.4.4]{PP03} the twisted structure sheaf
$\sO_{S_0}(\Theta_0)$ is $M$-regular, and  
$h^0 (S_0, \sO_{S_0}(\Theta_0) \otimes P_\xi) = 1$ for $P_\xi \in \Pic0(J(C_0))$ general. 
Thus the upper semicontinuity of the fibre dimension implies that for $t \neq 0$ sufficiently small, 
the $i$-th cohomological support locus
of $\sO_{S_t}(\Theta_t)$ has codimension strictly bigger than $i$ and 
$h^0 (S_t, \sO_{S_t}(\Theta_t) \otimes P_\xi) = 1$ for $P_\xi \in \Pic0(A_t)$ general. 
We have reached a contradiction. 
$\square$

\end{document}